\documentclass[ijoc]{informs3} 

\OneAndAHalfSpacedXII 
 \usepackage{booktabs}
 \usepackage{multirow}
 \usepackage{lscape}
 
\usepackage{natbib}
 \bibpunct[, ]{(}{)}{,}{a}{}{,}%
 \def\bibfont{\small}%
\TheoremsNumberedThrough     

\EquationsNumberedThrough 

\usepackage{graphicx,lmodern,amsfonts,mathrsfs,amsmath,amssymb,latexsym,xspace,epsfig,psfrag,color,tikz,tikzscale,pgfplots,colortbl}
\usepackage[flushleft]{threeparttable}
\usepackage[colorlinks,bookmarksnumbered,citecolor=blue]{hyperref}
\usepackage[ruled]{algorithm2e}
\usepackage[utf8]{inputenc}  
\usepackage{scalefnt}
\usepackage{subfigure}
\usepackage{multirow}

\usepackage{pdflscape}
\usepackage{multicol}

\usepackage{hyperref}
\usepackage{lineno}
\usepackage{longtable, lscape} 
\usepackage{bbm}
\usepackage{calrsfs}
\usepackage{txfonts}
\usepackage{rotating}
\usepackage{comment}

\usepackage{listings}
\definecolor{dkgreen}{rgb}{0,0.6,0}
\definecolor{gray}{rgb}{0.5,0.5,0.5}
\definecolor{mauve}{rgb}{0.58,0,0.82}

\usepackage{booktabs}
\usepackage{multirow}
\usepackage{pgfplotstable}
\usepackage{filecontents}
\usepackage{xspace}
\usepackage{tabularx}
\usepackage{siunitx}
\usepackage{lipsum}
\usepackage{tikz,pgfplots}
\usetikzlibrary{external}
\usetikzlibrary{patterns}
\pgfplotsset{compat=1.15}
\usepgfplotslibrary{statistics}
\usetikzlibrary{shapes,arrows}

\definecolor{Gray}{gray}{0.95}
\definecolor{Gray1}{gray}{0.6}
\definecolor{Gray2}{gray}{0.8}
\definecolor{Gray3}{gray}{1}

\begin{document}

\TITLE{AILS-II: An Adaptive Iterated Local Search Heuristic for the Large-scale Capacitated Vehicle Routing Problem}

\RUNAUTHOR{Máximo et al.}

\RUNTITLE{AILS-II: An Adaptive Iterated Local Search Heuristic for the Large-scale Capacitated Vehicle Routing Problem}

\ARTICLEAUTHORS{%
\AUTHOR{Vinícius R. Máximo}
\AFF{Instituto de Ciência e Tecnologia, Universidade Federal de São Paulo (UNIFESP)\\ Av. Cesare M. G. Lattes, 1201, Eugênio de Mello, São José dos Campos-SP, CEP: 12247-014, Brazil, \EMAIL{vinymax10@gmail.com}, \URL{}}
\AUTHOR{Jean-François Cordeau}
\AFF{HEC Montréal and GERAD, 3000 chemin de la Côte-Sainte-Catherine, Montréal, H3T 2A7, Canada, \EMAIL{jean-francois.cordeau@hec.ca}, \URL{}}
\AUTHOR{Mariá C. V. Nascimento}
\AFF{Divisão de Ciência da Computação (IEC), Instituto Tecnológico de Aeronáutica (ITA), Praça Marechal Eduardo Gomes, 50, Vila das Acácias, São José dos Campos-SP, CEP 12228-900, Brazil, \EMAIL{mariah@ita.br} \URL{} }
} 

\ABSTRACT{
A recent study on the classical Capacitated Vehicle Routing Problem (CVRP) introduced an adaptive version of the widely used Iterated Local Search (ILS) paradigm, hybridized with a path-relinking (PR) strategy. The solution method, called AILS-PR, outperformed existing meta-heuristics for the CVRP on benchmark instances. However, tests on large-scale instances suggest that PR is too slow, making AILS-PR less advantageous in this case. To overcome this challenge, this paper presents an Adaptive Iterated Local Search (AILS) combined with mechanisms to handle large CVRP instances, called AILS-II. The computational cost of this implementation is reduced while the algorithm also searches the solution space more efficiently.  AILS-II is very competitive on smaller instances, outperforming the other methods from the literature with respect to the average gap to the best known solutions. Moreover, AILS-II consistently outperforms  the state of the art on larger instances with up to 30,000 vertices. 
}



\maketitle
 


\section{Introduction}

The Capacitated Vehicle Routing Problem (CVRP) is a widely studied combinatorial optimization problem first introduced by \cite{Dantzig1959}. It consists in finding a set of routes that minimizes the cost of making deliveries to a set of customers with a homogeneous vehicle fleet based at a central depot. There is a great interest in providing efficient solutions to this problem, given its important role in many supply chains. \cite{Toth2001} carried out a study that reported potential savings of around 5 to 20\% in total transportation costs with the use of computer systems in real routing problems. Most of these problems are difficult to solve in reasonable time for large instances with hundreds or thousands of customers. In these cases, heuristic methods are often preferred as they provide  high-quality feasible solutions in reasonable computing time.

A number of successful heuristic methods have been employed to tackle the CVRP, among which we highlight the recently proposed Hybrid Genetic Search (HGS) \citep{Vidal2022}, Slack Induction by String Removals (SISRs) \citep{Christiaens2020}, Fast ILS Localized Optimization (FILO) \citep{Accorsi2021} and Adaptive Iterated Local Search with Path-Relinking (AILS-PR) \citep{Maximo2021}. In particular, \textcolor{black}{\cite{Maximo2021} demonstrated that AILS-PR outperformed  state-of-the-art metaheuristics published up to 2021} when considering small- and medium-sized instances, with up to 1000 vertices, proposed by \cite{Uchoa2017}, \cite{Christofides1979} and \cite{Golden1998}. However, recent experiments indicate that AILS-PR is too slow on  large-scale instances. More specifically, the path-relinking (PR) strategy shows poor performance in these cases. For large-scale CVRP instances, we highlight that FILO has achieved outstanding results.


This paper introduces a new version of AILS, called AILS-II, specially designed to efficiently solve large-scale CVRP instances. Several components were modified from the previous version of AILS to allow the processing of a large amount of data. Computational experiments with small- and medium-size instances show that AILS-II is very competitive with FILO and HGS. In addition, the experiments on large-scale instances with more than 336 customers indicate that AILS-II is consistently better than the other algorithms. AILS-II  is available  at the IJOC GitHub repository \citep{Maximoetal2024}.

The rest of the paper is organized as follows. Section \ref{sec:problem} presents the description of the CVRP and the main notations used throughout the paper. Section \ref{sec:proposedAILS} introduces AILS-II  and presents different approaches for controlling the perturbation degree and acceptance criteria. Section~\ref{codeorganization} shows a summary of the code organization and presents all parameters required to run the algorithm. Section~\ref{example} presents an example of how to solve a CVRP instance. Section~\ref{sec:EC} reports the results of computational experiments comparing AILS-II with other algorithms from the literature. Finally, in Section \ref{sec:conclusao}, we draw some conclusions and suggest future research directions.

\section{Problem Description}\label{sec:problem}

The Capacitated Vehicle Routing Problem (CVRP) can be defined on an undirected graph $G=(V,E)$, where $V = \{0,1,\ldots,n\}$ is the set of $n+1$ vertices and $E$ the set of edges. The depot is represented by the vertex $0$ and the other vertices in the subset $V_c = V\setminus \{0\}$ represent the $n$ customers. Each edge $\{i,j\} \in E$ has a non-negative weight $d_{ij}$ that represents the cost associated with a vehicle moving between vertices $i$ and $j$.

Each customer $ i \in V_c$ has a non-negative demand $q_i$ that must be satisfied (the depot has demand $q_0=0$). To meet the demands of the customers, $m$ identical vehicles with a capacity of $\bar{q}$ are used. To ensure the feasibility of the problem, the demand of each customer is assumed to be smaller than or equal to the vehicle capacity, that is, $q_i \leq \bar{q} ~\forall i \in V_c$. In the CVRP, each route takes the form of a closed loop with no node repetition. A closed loop is represented by a cyclic sequence of vertices where a pair of vertices is adjacent if they are consecutive in the sequence and not adjacent otherwise. The CVRP, therefore, consists in finding a set of $m$ routes that minimize the sum of the edge weights. The routes of a solution $s$ are described by $\mathscr{R}=\{R^s_1,R^s_2,\ldots, R^s_{m_s}\}$, where $m_s$ represents the number of routes in the solution $s$, $R^s_i = \{v^{i}_0, v^{i}_1,\ldots,v^{i}_{t_i}\}$, $t_i + 1$ is the length of route $R^s_i$, $v^{i}_0= v^{i}_{t_i}=0$, $R^s_i \cap R^s_j =\{0\}$, for $i\neq j$, and $\cup_{i=1}^{m_s} R^ s_i = V$. 

The CVRP is $\mathbb{NP}$-Hard \citep{Lenstra1981} since it generalizes the Traveling Salesman Problem (TSP), which seeks to minimize the length of a Hamiltonian tour. 

\section{Adaptive Iterated Local Search (AILS) \textcolor{black}{For Large-Scale Instances}}
\label{sec:proposedAILS}

AILS was initially proposed by \cite{Maximo2021} to solve the CVRP, and has recently been adapted for the Heterogeneous Fleet Vehicle Routing Problem (HVRP) by \cite{Maximo2022}. AILS is an adaptive metaheuristic that has two main steps: perturbation and local search. These two steps are performed iteratively until a stopping criterion is reached. At each iteration, a reference solution is perturbed to generate a potentially different solution. A local search step is then applied to improve the quality of the resulting solution through an exploration of the neighborhood formed by vertex and edge movements. Then, this solution is evaluated by the acceptance criterion. If it is accepted, it becomes the new reference solution. The solution obtained after the perturbation step may be infeasible. For this reason, we apply an algorithm that uses the same neighborhood as the local search but will guarantee the feasibility of the solution. The main difference between the local search and feasibility strategies is that the movements in the latter are chosen by prioritizing the reduction of capacity constraint violations. The adaptive behavior of AILS can be observed in the definition of the perturbation degree and in the acceptance criterion, described in more detail in Sections \ref{ControlePert} and \ref{criterioAceitacao}, respectively. These two conditions are responsible for the diversity control of the method, which means that controlling them is of utmost importance. \textcolor{black}{A general pseudocode of AILS-II is presented in Algorithm~\ref{alg:AILS}.}


\begin{algorithm}[H]
\LinesNumbered
\SetAlgoLined
\KwData{Instance data}
\KwResult{The best solution found $s^*$}
$s \leftarrow$ Construct an initial solution \\ 
$s^r,s^* \leftarrow$ Local Search($s$) \\
\Repeat{stopping criterion is met}
{ 
    $s \leftarrow$ Perturbation Procedure  ($s^r$)\\
    $s \leftarrow$ Local Search($s$)\label{AILSPR:linha5ails}\\
	Update the diversity control parameters considering the distance between $s$ and $s^r$\\
    $s^r\leftarrow $ Apply acceptance criterion to $s$\\
Update the acceptance criterion\\
	Assign $s$ to $s^*$ if $ f(s) < f(s^*)$\\

}

\caption{Adaptive Iterated Local Search}
\label{alg:AILS}

\end{algorithm}

\subsection{Perturbation Degree}
\label{ControlePert} 

The perturbation degree establishes the number of changes applied to the reference solution $s^r$ to obtain a different solution. The greater the perturbation degree, the larger the distance between these solutions. Thus, the control of the perturbation degree enables the method to manage the diversity of the search algorithm. For this reason, the perturbation degree control can be seen as a mechanism of great relevance for metaheuristics, considering that adequate control can allow the algorithm to escape from local optima. High diversity means that the algorithm will be able to escape from a local optimum more easily, but it might be costly to find the best solution in a given neighborhood \citep{Lourenco2010}. On the other hand, low diversity yields a higher chance that the algorithm will get stuck in a local optimum. To achieve an adequate balance between these two goals,  AILS-II considers adaptive mechanisms to adjust such a parameter. The mechanism proposed in \citet{Maximo2021} uses a fixed parameter, called $d_{i}$, that establishes the ideal distance between the solution $s$ obtained by the local search and the reference solution $s^r$. For AILS-II, we introduce a convergent behavior to control the degree of perturbation.  In this convergent mechanism, the value of $d_{i}$ starts with the value of $d_{max}$ and decreases throughout the execution of the algorithm until it reaches a value of $d_{min}$. For this, at each iteration of the algorithm, the value of $d_{i}$ is adjusted as $d_{i} \leftarrow d_{i} \left( \frac{d_{min}}{d_{max}}\right)^{\frac{1}{it}} $, where $it$ is the estimated maximum number of iterations performed by the algorithm.

\subsection{Acceptance Criterion}
\label{criterioAceitacao}

The acceptance criterion establishes whether the current solution should become the reference solution and be used in the following iterations of the algorithm. The AILS introduced in this paper uses a convergent acceptance criterion with a more relaxed criterion at the beginning of the search and a more restrictive condition as it approaches the end of the search. This criterion was inspired by the Threshold Acceptance (TA) algorithm proposed by \cite{Dueck1990}. In line with this, the employed acceptance criterion  restricts the quality of the accepted solutions to the threshold $\theta = \underline{f} + \eta (\bar{f}-\underline{f})$. \textcolor{black}{This means that $\theta$ can assume values in the interval [$\underline {f}$,$\bar{f}$], where $\underline {f}$ is the best solution found in the last $\gamma$ iterations, and $\bar{f}$ is the average quality of the solutions obtained by the local search. The value of $\eta$ establishes how far $\theta$ is from the lower and upper limits of the interval. If $\eta$ takes a value close to 1, then $\theta$ is more relaxed, with values closer to $\bar{f}$. Otherwise, values close to 0 for $\eta$ make $\theta$ closer to $\underline {f}$.} This threshold was proposed by \citet{Maximo2021} and the $\eta$ value was adjusted according to the flow of accepted solutions. In AILS-II, we propose a convergent variation for the value of $\eta$. Therefore, $\eta$ starts at 1 and decreases  until a minimum value $\epsilon=0.01$. Thus, at each iteration of the algorithm $\eta \leftarrow \eta \epsilon ^{\frac{1}{it}}$.

\subsection{\textcolor{black}{More Particularities of AILS-II}}

Besides the convergent criteria in the diversity control mechanisms, this version of AILS presents substantial changes from its previous versions. These were necessary to achieve a good performance on large instances. We enumerate next the key differences between AILS-II and previous versions of AILS:
\begin{itemize}
    \item \textbf{\textcolor{black}{New combination of movements:}}  \textcolor{black}{The configuration of the local search and feasibility phases has the intra-route movements SHIFT, SWAP, 2-opt and the inter-route movements SHIFT, SWAP* \citep{Vidal2022} and 2-opt*. Therefore, a difference from the AILS-PR \citep{Maximo2021} is the replacement of SWAP by SWAP* in the inter-route movement. Moreover, to achieve a better efficiency, we restricted the neighborhood in the SWAP* by considering for each vertex its $\varphi$ closest vertices and not the entire set of vertices that belong to its route. The motivation is that the larger the route size, the greater the computational cost. For large instances, the computational cost becomes prohibitive.}

\item \textbf{Feasibility algorithm:} AILS-PR uses two different criteria to assign values to rank the moves depending on whether the performed move will provide a better or worse quality solution. AILS-II  considers the same criterion for both cases and ranks the moves according to the ratio of the difference in cost between the original solution and the solution after the move and the feasibility gain. The feasibility gain is the value of the infeasibility reduction considering the violated constraints.

\item \textcolor{black}{\textbf{Limited neighborhood of the local search and feasibility procedures}: To speed up the local search and the feasibility procedures, we imposed in AILS-II that their movements are applied only in modified routes, i.e., routes that were perturbed. In the case of the inter-route movements,  one of the routes must necessarily be the one that has been perturbed.}

\item \textbf{Different perturbation heuristics:} The perturbation heuristics of AILS are composed of a set of addition and removal heuristics. Removal heuristics remove vertices from a given solution whereas addition heuristics add them back in positions that depend on the criterion of the method. Both AILS-PR and AILS-II follow the batch approach to remove and add them back to the solutions. This means that they first remove $\omega_k$ vertices from the solution through a removal heuristic $R_k$ and, after that, the addition heuristic inserts the vertices back in the solution. The AILS for the HVRP follows an alternate call between removal and addition heuristics, where one of the $\omega_k$ vertices is removed and subsequently added to a position of the solution. The process repeats $\omega_k$ times, to allow all $\omega_k$ vertices to change position. On the one hand, AILS-II does allow the vertices to return to the same position. On the other hand, with the addition heuristics  in AILS-PR and AILS for the HVRP, vertices can be inserted in the same position. Regarding the employed removal heuristics, all of them use the concentric and sequence removal strategies. However, AILS-PR considers two other removal heuristics, whereas AILS for the HVRP also employs the random removal. AILS-II considers two addition heuristics. The only addition heuristic common to AILS for the HVRP, AILS-PR  and AILS-II is the one called insertion by cost. Besides, AILS-II uses the insertion by distance, also considered by AILS for the HVRP. \textcolor{black}{We refer to \cite{Maximo2021} and \cite{Maximo2022} for a better description of the perturbations approaches and methods.} Table \ref{versoes}  presents the main heuristics used in the perturbation step of the three different AILS versions as well as the employed approaches. 
\begin{table}[!ht]
\caption{Main characteristics of the perturbation step in the AILS-PR \citep{Maximo2021}, AILS for the HVRP \citep{Maximo2022} and AILS-II.}
\label{versoes}
\centering
\begin{tabular}{lcccc}
\hline
\hline
 & & \multicolumn{3}{c}{Algorithm} \\ \cline{2-5}
Perturbation approach/method    &&  AILS-PR (CVRP)  & AILS  (HVRP)  & AILS-II (This paper)     \\
\hline
\hline
\rowcolor{Gray}
Batch approach &&       \checkmark &               & \checkmark  \\
Alternate approach&&               & \checkmark    &             \\
\rowcolor{Gray}
Concentric removal && \checkmark    & \checkmark    & \checkmark  \\
Random removal &&                & \checkmark    &             \\
\rowcolor{Gray}
Sequence removal &&  \checkmark & \checkmark    & \checkmark  \\
Proximity removal &&  \checkmark    &               &             \\
\rowcolor{Gray}
Insertion by proximity && \checkmark&               &             \\ 
Insertion by cost &&     \checkmark & \checkmark    & \checkmark  \\
\rowcolor{Gray}
Insertion by distance &&            & \checkmark    & \checkmark  \\
\hline
\hline
\end{tabular}
\end{table}



\end{itemize}
  

  

 AILS-II does not employ the so-called path-relinking (PR), present in AILS-PR \citep{Maximo2021}. The reason for it is that the computational cost required by the hybridized version on large-scale CVRP instances is too high. 

\begin{figure}
    \centering
\includegraphics[]{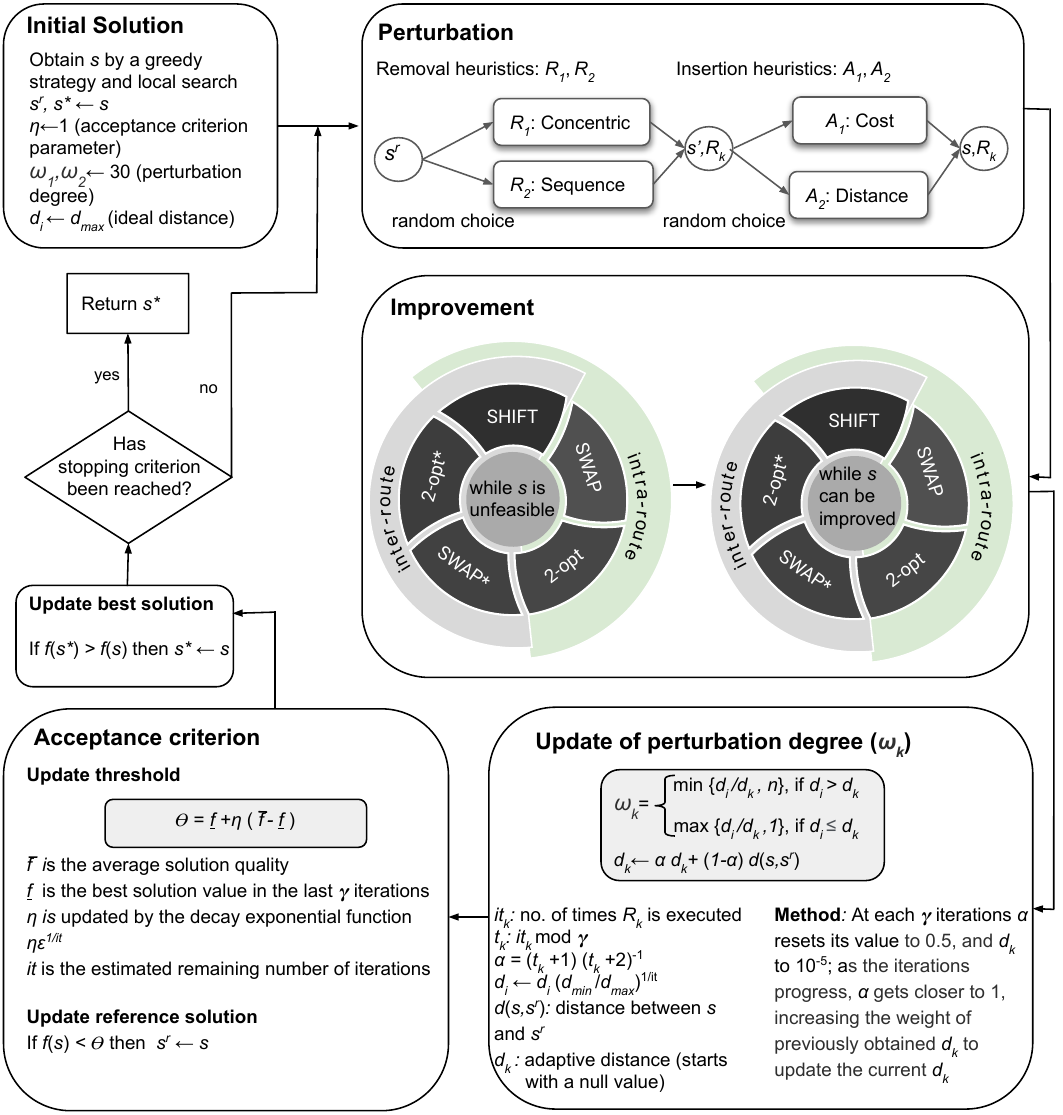}
    \caption{Flowchart of AILS-II}
    \label{fig:flowchartAILS}
\end{figure}

\textcolor{black}{A diagram summarizing all the steps of AILS-II is presented in Figure~\ref{fig:flowchartAILS}.}


\section{Code Organization}\label{codeorganization}

AILS-II was developed in Java and the code is available at \url{https://github.com/INFORMSJoC/2023.0106}. The code is divided into 8 packages containing 35 files and roughly 6000 lines. The main class is called AILSII, which is inside the \texttt{SearchMethod} package. The main method is presented in Source Code  \ref{main}. The user parameter values are input provided by the user. 
An instance of the problem is constructed and used to create the \texttt{ailsII} object from the \texttt{AILSII} class. This class receives in the constructor an instance of the problem and the configuration of the algorithm. The beginning of the search occurs with the call to the search method of the \texttt{AILSII}  class.

\begin{lstlisting}[
language=Java, 
caption=Main method of AILS,
label=main
]
public static void main(String[] args) 
{
	InputParameters reader=new InputParameters();
	reader.readingInput(args);
	Instance instance=new Instance(reader);
	AILSII ailsII=new AILSII(instance, reader);
	ailsII.search();
}
\end{lstlisting}


    

A summary of the contents of each package is presented next.

\begin{itemize}
    \item \texttt{Auxiliary}: This package contains the method that calculates the distance between two solutions and the dynamic average used in the diversity control.
    \item \texttt{Evaluators}: This package contains the cost and feasibility evaluators of the movements, and the execution process of the movements used by the local search and feasibility.
    \item  \texttt{DiversityControl}: This package contains the source codes of the acceptance criterion and the control of the perturbation degree.
    \item \texttt{Data}: This package contains codes for reading instance data.
    \item \texttt{Improvement}: This package contains the implementation of the local search and feasibility methods.
    \item \texttt{SearchMethod}: This package contains the implementation of AILS-II.
    \item \texttt{Perturbation}: This package contains  the addition and removal heuristics  used in the perturbation step.
    \item \texttt{Solution}: This package contains the data structures of the solution, routes and graph vertices.
    
\end{itemize}

\section{Example of Usage} \label{example}


To run the AILSII class it is necessary to define the following parameters:

\begin{itemize}
    \item \texttt{-file}: the file name of the problem instance.
 
    \item \texttt{-rounded}: A flag that indicates whether the instance has rounded distances or not. The options are: [false, true]. The default value is true.
    
    \item \texttt{-stoppingCriterion}: It is possible to use two different stopping criteria:
    \begin{itemize}
        \item \texttt{Time}: The algorithm stops when a given time in seconds has elapsed.
        \item \texttt{Iteration}: The algorithm stops when the given number of iterations has been reached. 
    \end{itemize}

    \item \texttt{-limit}: Refers to the value that will be used in the stopping criterion. If the stopping criterion is a time limit, this parameter is the timeout in seconds. Otherwise, this parameter indicates the number of iterations. The default value is the maximum limit for a double precision number in the JAVA language.

    \item \texttt{-best}: Indicates the value of the optimal solution. The default value is 0.

    \item \texttt{-varphi}: Parameter of the feasibility and local search methods that refers to the maximum cardinality of the set of nearest neighbors of the vertices. The default value is 40. The larger it is, the greater the number of movements under consideration in the methods. 

    \item \texttt{-gamma}: Number of iterations for AILS-II to perform a new adjustment of $\omega$. The default value is 30.

    \item \texttt{-dMax}: Initial reference distance between the reference solution and the  solution obtained after the local search. The default value is 30.

    \item \texttt{-dMin}: Final reference distance between the reference solution and the solution obtained after the local search. The default value is 15.

 \end{itemize}

An example of how to execute AILS-II would be:

\begin{verbatim}
java -jar AILSII.jar -file Instances/X-n214-k11.vrp -rounded true 
-best 10856 -limit 100 -stoppingCriterion Time 

\end{verbatim}
This command will run AILS-II to solve an instance whose parameters are in file ``Instances/X-n214-k11.vrp", considering that the distance between the vertices is rounded to an integer value. We also indicate that the optimal value for this instance is 10856, therefore, if the algorithm finds a solution with this cost, it halts. The stopping criterion is a time limit of 100 seconds.

\section{Computational Experiments}
\label{sec:EC}

The computational experiments were performed on a cluster with 104 nodes, each of them with 2 Intel Xeon E5-2680v2 processors running at 2.8 GHz, 10 cores, and 128 GB DDR3 of 1866 MHz RAM. In the first experiment, we investigate the performance of the proposed algorithm using the set of benchmark instances presented in \cite{Uchoa2017}. This set contains 100 CVRP instances whose size ranges from 100 to 1000 vertices. The second experiment assesses the performance of AILS-II on larger instances. We consider the 10 instances proposed in \cite{Arnold2019} which have between 3000 and 30,000 vertices. All experiments were carried out using as stopping criterion $3n$ seconds, where $n$ is the number of vertices in the instance.

\subsection{First Experiment}

This experiment contrasts the performance of AILS-II to HGS \citep{Vidal2022} and FILO \citep{Accorsi2021} on the 100 instances proposed by \cite{Uchoa2017}. We run 50 times each of them and report the average gap, defined as:

\begin{equation}
   gap=100\times \frac{(Avg-\mbox{BKS})}{\mbox{BKS}},\label{avggap}
\end{equation} 

\noindent where BKS is the objective function value of the best known solution of a given instance and $Avg$ is the average objective function value of the solutions obtained in the independent executions.

\begin{figure}[htbp]
\center
\includegraphics[]{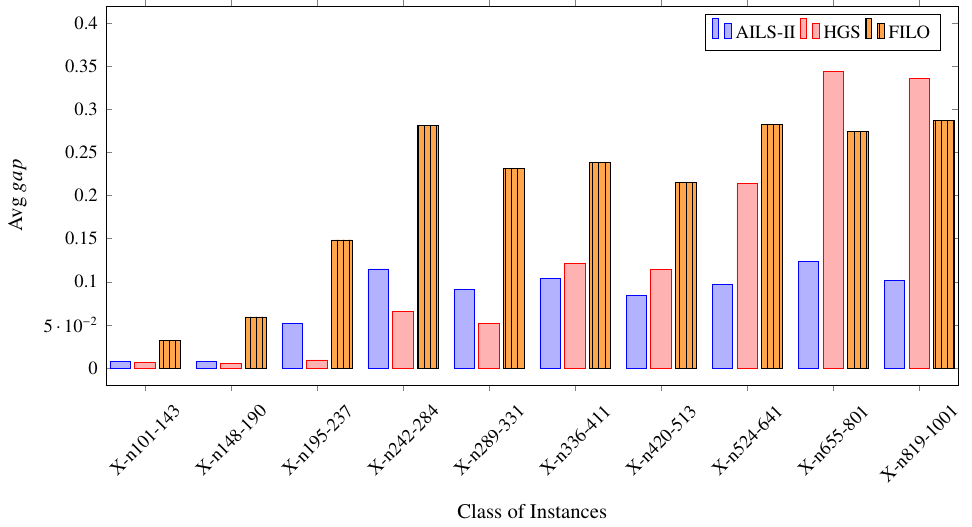}
\caption{Comparison of the average gaps obtained by the AILS-II, HGS and FILO  on the 100 instances proposed by \cite{Uchoa2017}.}
\label{fig:gap}
\end{figure}

Figure \ref{fig:gap} shows the average gap of 10 classes of instances that were grouped by size. Considering instances with up to 331 vertices, we observed the following order of best performance: HGS, AILS-II and FILO. In instances with more than 336 vertices, AILS-II achieves the best results and for instances with more than 655 vertices, FILO has a better result than HGS. AILS-II presented a lower computational cost. In these experiments, AILS-II found a better best known solution for the ten following instances: X-n536-k96 = 94828, X-n701-k44=81919, X-n716-k35=43330, X-n766-k71=114416, X-n783-k48=72381, X-n837-k142=193734, X-n895-k37=53848, X-n916-k207=329178, X-n957-k87=85464, X-n979-k58=118913. \textcolor{black}{Appendix~\ref{appendix} displays the numerical results for each of the instances used in this experiment.}

\begin{figure}[htbp]
\center
\includegraphics[]{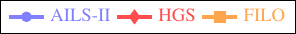}
\label{legenda}
\\ 
\includegraphics[]{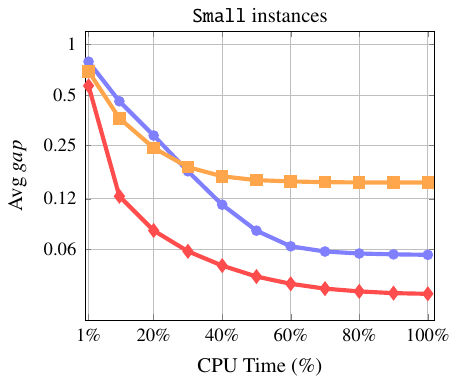}
\includegraphics[]{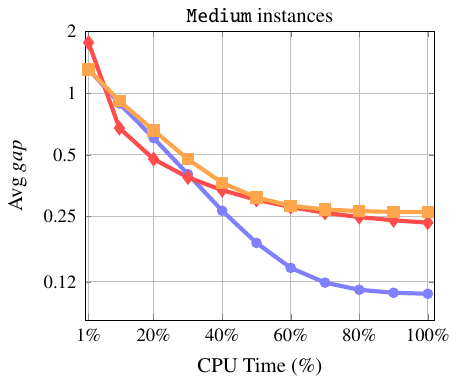}
\includegraphics[]{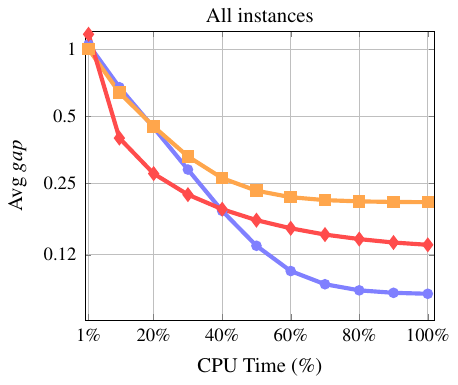}
\caption{\textcolor{black}{Convergence of the  AILS-II, HGS and FILO over time for the X instances \citep{Uchoa2017}. The plots refer to, respectively, the average  results on the 50 smallest instances (with up to 331 nodes), on the 50 largest instances (with more than 335 nodes) and on all instances.}}
\label{fig:convergence}
\end{figure}

\textcolor{black}{An experiment showing the convergence over time of all tested algorithms was performed for a better analysis of AILS-II regarding existing algorithms. In line with this, as in \cite{Vidal2022}, we divided the X instances into two groups. The first is referred to as \texttt{Small} and contains the first 50 instances, with 100 to 330 customers. The second group is called \texttt{Medium}, with the second half of X instances, which include from 335 to 1000 customers. Figure~\ref{fig:convergence} displays the average gap (Avg \textit{gap}) when considering 10\%, 20\%, 30\%, 40\%, 50\%, 60\%, 70\%, 80\%, 90\% and 100\% of the time imposed for all algorithms. Besides the average results for \texttt{Small} and \texttt{Medium} instances, we also present results for all instances. These plots show a clear behavior of the algorithms over time on the two groups of instances. AILS-II was the second best on \texttt{Small} instances when considering 30\% to 100\% of the CPU time. Despite being very competitive with FILO when considering a CPU time up to 30\%, FILO and HGS were better than AILS-II. On the other hand, for \texttt{Medium} instances, AILS-II achieved the second best results when imposing 5\% to 30\% of the CPU time, but with a slight difference from FILO. However, when considering over 50\% of the time, AILS-II outperformed the other algorithms achieving half of their average gap in 100\% of the CPU time. Regarding the average gap on all instances, AILS-II required over 40\% of the CPU time to outperform the other algorithms.}
 
\subsection{Second Experiment}

In this experiment, 10 independent executions of AILS-II and FILO were performed on the instances proposed by \cite{Arnold2019}. HGS was not considered in this experiment because its author suggested not to use it \citep{Vidal2022}. The reason is that these instances have very different characteristics from those for which the algorithm was designed. 

Table~\ref{Arnold} shows, for each instance, the instance name, the number of nodes ($n$),  the best-known solution, followed by the results achieved by FILO and AILS-II. The results compiled in the table are the average solution value of the independent executions (\textit{Avg}), the average gap ($gap$), the  objective function value of the best solution obtained in the runs (Best); and the average time in minutes that the algorithm took to find the best solution in each run (T (min)). BKSs highlighted with `*' were found by AILS-II and improve previous BKSs. 
It is worth mentioning that there were some BKSs identified when setting up the parameter values of AILS-II. That is the reason why some BKSs improved by AILS-II (not necessarily in the final experiment) were better than the best solution achieved by AILS-II in the runs with the final parameter values we reported in this table.
\begin{table}[!ht]
\centering
\scalefont{0.75}
\caption{Results achieved by AILS-II and FILO on the instances introduced by \cite{Arnold2019}.}
\label{Arnold}
\centering
\begin{tabular}{lcccccccccc}
\hline
&&& \multicolumn{3}{c}{FILO \citep{Accorsi2021}} && \multicolumn{3}{c}{AILS-II}\\ 
\cline{4-6} \cline{8-10}  
Instance & $n$ & BKS&\textit{Avg}~($gap$)&Best&T (min)&&\textit{Avg}~($gap$)&Best&T (min)\\
\hline
\rowcolor{Gray}
Antwerp1&6000&$477277$&478002.8~(0.1521)&477854&298.3&&\textbf{477526.5~(0.0523)}&\textbf{477466}&\textbf{289.0}\\
Antwerp2&7000&$291350$&291776.6~(0.1464)&\textbf{291511}&348.0&&\textbf{291715.8~(0.1256)}&291587&\textbf{338.4}\\
\rowcolor{Gray}
Brussels1&15000&$501584^*$&502349.1~(0.1525)&502175&746.5&&\textbf{501883.5~(0.0597)}&\textbf{501735}&\textbf{737.8}\\
Brussels2&16000&$345057^*$&345940.4~(0.2560)&345691&798.3&&\textbf{345419.0~(0.1049)}&\textbf{345253}&\textbf{792.3}\\
\rowcolor{Gray}
Flanders1&20000&$7238970^*$&7249408.7~(0.1442)&7248487&996.3&&\textbf{7240768.8~(0.0248)}&\textbf{7239443}&\textbf{982.6}\\
Flanders2&30000&$4367291^*$&4384803.7~(0.4010)&4382783&1494.8&&\textbf{4370437.6~(0.0720)}&\textbf{4369014}&\textbf{1477.3}\\
\rowcolor{Gray}
Ghent1&10000&$469483^*$&469991.9~(0.1084)&469830&496.5&&\textbf{469684.0~(0.0428)}&\textbf{469476}&\textbf{488.5}\\
Ghent2&11000&$257563^*$&258173.0~(0.2368)&258052&546.9&&\textbf{257926.1~(0.1410)}&\textbf{257670}&\textbf{543.1}\\
\rowcolor{Gray}
Leuven1&3000&$192848$&193119.1~(0.1406)&193014&146.0&&\textbf{193013.1~(0.0856)}&\textbf{192923}&\textbf{143.4}\\
Leuven2&4000&$111391$&111706.2~(0.2830)&111618&197.7&&\textbf{111703.0~(0.2801)}&\textbf{111512}&\textbf{192.0}\\
\rowcolor{Gray}
Average&&&~(0.2021)&&606.9&&~(0.0989)&&598.4\\
\hline
\end{tabular}
\end{table}
\addtolength{\tabcolsep}{3pt}

AILS-II achieved the smallest average gap in all 10 instances. For one instance we did not find the best solution. These results confirm the better performance of AILS-II on large instances with respect to the state-of-the-art large scale heuristic method.

\begin{figure}[htbp]
\center
\pgfplotsset{footnotesize,samples=10}
\includegraphics[]{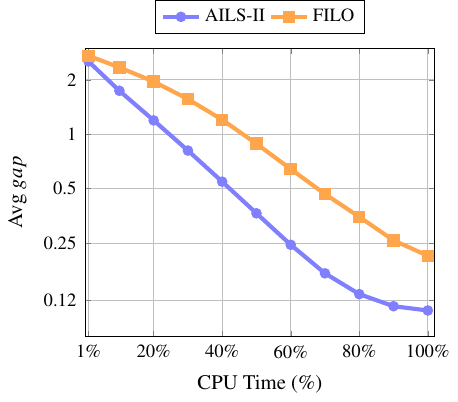}
\caption{\textcolor{black}{Convergence of the AILS-II and FILO over time for the large-scale instances. This plot refers to the average results on the instances introduced by \cite{Arnold2019}.}}
\label{fig:convergence_large}
\end{figure}

\textcolor{black}{Similar to the convergence plots presented in the previous experiment, Figure~\ref{fig:convergence_large} displays the average gap when considering from 10\% to 100\% of the time with a step of 10\%. The plot shows that AILS-II achieved the best gaps within all time limits, demonstrating its high performance on large-scale instances.}

\subsection{\textcolor{black}{Performance of the main components of AILS-II}}

\begin{figure}[h!t]
\center

\includegraphics[]{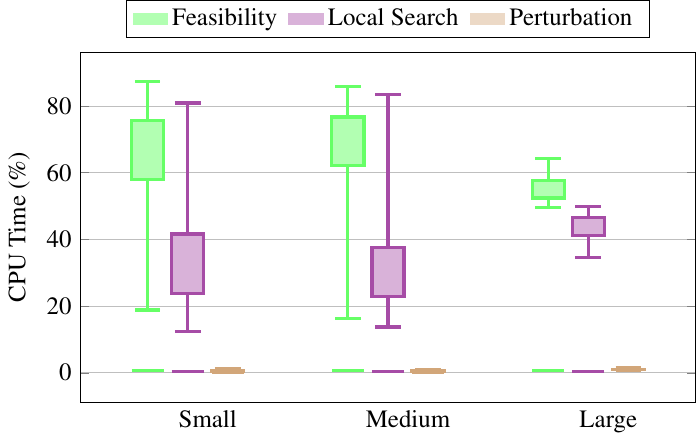}
\caption{Box plots of the percentage of time spent by the different parts of AILS-II.}
\label{fig:boxplot}
\end{figure}

\textcolor{black}{This section presents an experiment displaying the statistics of the time required in each of the primary components of AILS-II: perturbation, feasibility and local search. Figure~\ref{fig:boxplot} shows box plots of the amount of time required when testing X instances, divided on \texttt{Small} and \texttt{Medium}, and on the large-scale instances.}

\textcolor{black}{In all tests, the most time-consuming part of AILS-II is the feasibility phase, which requires over 50\% of the execution time  regardless of the type of tested instances. One can also observe that the percentage of the time spent in the different parts of AILS-II on the X instances were very similar for small and medium size instances. On large-scale instances, the time required for the local search and the feasibility phases was more balanced. In all cases, as expected, the time for perturbation is marginal in comparison to the other two phases of AILS-II.}

\textcolor{black}{Finally, in Appendix~\ref{appendix2} we also report the results of the sensitivity analysis on the four main parameters controlling the heuristic.}

\section{Final Remarks and Future Work}
\label{sec:conclusao}
 
AILS is an adaptive ILS metaheuristic introduced  by \cite{Maximo2021} to solve the CVRP. The authors proposed a hybridization of the path-relinking strategy with AILS, called  AILS-PR, to better explore the solution space of such a problem. As a result, AILS-PR achieved outstanding results on the tested instances, which consisted mostly of small- and medium-size instances. However, preliminary experiments indicated that neither AILS-PR nor such a version of AILS without PR performed well on large-scale instances. On the one hand, AILS does not intensify the search much. On the other hand, AILS-PR is too slow on large instances.

Therefore, the main goal of this study was to present a version of AILS, called AILS-II, more efficient approach for large instances of the CVRP. \textcolor{black}{For this,  we were very thorough in implementing the different components of the metaheuristic. The local search, which is one of the most time-consuming parts of the method, is composed of movements  common to the other state-of-the-art algorithms. Therefore, our effort to achieve an efficient method was directed to (i) the adequate fine tuning of the main ingredients of ILS: the perturbation degree and the acceptance criterion; and (ii) developing efficient data structures and designing low-cost strategies. This fine tuning considered a combination of over a thousand forms in order for the AILS-II to make the right decisions to jump from a solution to another,  avoid redundant processes,  and decide which are the most promising solutions to intensify the search (to become the reference solution).  AILS-II} employed a convergent criterion for both the diversity control strategy and in the acceptance criterion. The results show that the adequate control of these two stages in the adaptive process is of great relevance to ensure good performance.

The computational experiments were performed using a benchmark dataset proposed by \cite{Uchoa2017} that contains 100 instances whose size ranges from 100 to 1000 vertices. We compared AILS-II with the HGS and FILO. Considering this dataset, AILS-II achieved the best average gap in all instances with more than 336 vertices. To evaluate the performance of AILS-II for larger instances, we carried out experiments with a dataset proposed by \cite{Arnold2019} that contains 10 large instances of CVRP. These instances have between 3000 and 30,000 vertices. The results of this experiment showed that AILS-II outperformed FILO on all instances.

As future work, we intend to continue investigating more efficient mechanisms for controlling diversity. We also want to explore other search strategies that are robust and simpler. Other areas of interest involve machine learning to optimize the performance of metaheuristics. \textcolor{black}{Given that the code is fast and flexible (it was adapted to handle the heterogeneous fleet VRP by \cite{Maximo2022} and the maximal covering location problem by \cite{Maximo2023}), we recommend its use to researchers wishing to solve large-scale instances of the classical VRP or to adapt it to related problem variants.}

\section*{Acknowledgments}

\textcolor{black}{We are thankful to two anonymous referees for their valuable comments on an earlier version of this paper.} The authors are \textcolor{black}{also} grateful for the financial support provided by CNPq (403735/2021-1; 309385/2021-0) and FAPESP (2019/22067-6; 2022/05803-3); and the Center for Mathematical Sciences Applied to Industry (CeMEAI), funded by FAPESP by through the 2013/07375-0 project, for the availability of computer resources.

\bibliographystyle{informs2014} \bibliography{Bibliography}





\newpage

\begin{APPENDIX}{Other results}\label{appendix1}
\textcolor{black}{This appendix presents the numerical results on X instances and a brief parameter analysis of AILS-II.}
\section{Numerical results}\label{appendix}
Table~\ref{tab:xinstances} presents the results of FILO, HGS and AILS-II for the benchmark instances from \citep{Uchoa2017}. Column `BKS' indicates best known solutions highlighted with an asterisk when their values were improved in experiments performed in this paper. Columns `\textit{Avg}~($gap$)', `Best' and `T (min)' display, respectively, the average solution and mean gap between parentheses, the best solution and average time of the independent executions. 


\begin{landscape}
\footnotesize
\begin{longtable}{@{}lllllllllll@{}}
\caption{Results achieved by FILO \citep{Accorsi2021}, HGS \citep{Vidal2022} and AILS-II on the benchmark instances introduced by \citep{Uchoa2017}.}
\label{tab:xinstances}\\
\toprule
\multirow{2}{*}{Instance} & \multirow{2}{*}{BKS} & \multicolumn{3}{c}{FILO \citep{Accorsi2021}} & \multicolumn{3}{c}{HGS \citep{Vidal2022}} & \multicolumn{3}{c}{AILS-II} \\* \cmidrule(r){3-5}\cmidrule(r){6-8}\cmidrule(r){9-11}
 & & \textit{Avg}~($gap$) & Best & T (min) & \textit{Avg}~($gap$) & Best & T (min) & \textit{Avg}~($gap$) & Best & T (min) \\* \cmidrule(r){1-11}
\endfirsthead
\endhead
\bottomrule
\endfoot
\endlastfoot

X-n101-k25 & \textbf{$27591^*$} & \textbf{27591~(0.0000)} & \textbf{27591} & \textbf{0.02} & \textbf{27591.00~(0.0000)} & \textbf{27591} & \textbf{0.02} & \textbf{27591.00~(0.0000)} & \textbf{27591} & 1.01 \\
X-n106-k14 & \textbf{$26362^*$} & 26375.76~(0.0522) & \textbf{26362} & \textbf{1.13} & 26379.32~(0.0657) & \textbf{26362} & 2.70 & \textbf{26362.16~(0.0006)} & \textbf{26362} & 2.55 \\
X-n110-k13 & \textbf{$14971^*$} & \textbf{14971~(0.0000)} & \textbf{14971} & \textbf{0.01} & \textbf{14971.00~(0.0000)} & \textbf{14971} & \textbf{0.01} & \textbf{14971.00~(0.0000)} & \textbf{14971} & 0.22 \\
X-n115-k10 & \textbf{$12747^*$} & \textbf{12747~(0.0000)} & \textbf{12747} & \textbf{0.03} & \textbf{12747.00~(0.0000)} & \textbf{12747} & \textbf{0.03} & \textbf{12747.00~(0.0000)} & \textbf{12747} & 0.27 \\
X-n120-k6 & \textbf{${13332^*}$} & \textbf{13332~(0.0000)} & \textbf{13332} & 0.19 & \textbf{13332.00~(0.0000)} & \textbf{13332} & \textbf{0.16} & \textbf{13332.00~(0.0000)} & \textbf{13332} & 1.01 \\
X-n125-k30 & \textbf{${55539^*}$} & 55542.74~(0.0067) & \textbf{55539} & 0.56 & \textbf{55539.00~(0.0000)} & \textbf{55539} & \textbf{0.44} & \textbf{55539.00~(0.0000)} & \textbf{55539} & 3.49 \\
X-n129-k18 & \textbf{${28940^*}$} & 28949.94~(0.0343) & \textbf{28940} & 1.12 & \textbf{28940.00~(0.0000)} & \textbf{28940} & \textbf{0.16} & 28942.08~(0.0072) & \textbf{28940} & 2.20 \\
X-n134-k13 & \textbf{${10916^*}$} & 10924.3~(0.0760) & \textbf{10916} & 0.85 & \textbf{10916.00~(0.0000)} & \textbf{10916} & \textbf{0.53} & 10916.04~(0.0004) & \textbf{10916} & 1.73 \\
X-n139-k10 & \textbf{${13590^*}$} & \textbf{13590~(0.0000)} & \textbf{13590} & 0.55 & \textbf{13590.00~(0.0000)} & \textbf{13590} & \textbf{0.03} & 13590.08~(0.0006) & \textbf{13590} & 2.17 \\
X-n143-k7 & \textbf{${15700^*}$} & 15723.84~(0.1518) & \textbf{15700} & 1.93 & \textbf{15700.00~(0.0000)} & \textbf{15700} & \textbf{0.21} & 15710.84~(0.0690) & \textbf{15700} & 3.03 \\
X-n148-k46 & \textbf{${43448^*}$} & 43456.52~(0.0196) & \textbf{43448} & 1.08 & \textbf{43448.00~(0.0000)} & \textbf{43448} & \textbf{0.10} & \textbf{43448.00~(0.0000)} & \textbf{43448} & 2.38 \\
X-n153-k22 & \textbf{${21220^*}$} & 21237.5~(0.0825) & 21225 & 1.70 & 21224.88~(0.0230) & \textbf{21220} & \textbf{0.52} & \textbf{21224.62~(0.0218)} & \textbf{21220} & 3.16 \\
X-n157-k13 & \textbf{${16876^*}$} & \textbf{16876~(0.0000)} & \textbf{16876} & 1.68 & \textbf{16876.00~(0.0000)} & \textbf{16876} & \textbf{0.10} & \textbf{16876.00~(0.0000)} & \textbf{16876} & 2.75 \\
X-n162-k11 & \textbf{${14138^*}$} & 14156.62~(0.1317) & \textbf{14138} & 2.40 & \textbf{14138.00~(0.0000)} & \textbf{14138} & \textbf{0.05} & \textbf{14138.00~(0.0000)} & \textbf{14138} & 1.79 \\
X-n167-k10 & \textbf{${20557^*}$} & 20557.86~(0.0042) & \textbf{20557} & 1.22 & \textbf{20557.00~(0.0000)} & \textbf{20557} & \textbf{0.36} & 20559.92~(0.0142) & \textbf{20557} & 2.83 \\
X-n172-k51 & \textbf{${45607^*}$} & \textbf{45607~(0.0000)} & \textbf{45607} & 2.56 & \textbf{45607.00~(0.0000)} & \textbf{45607} & \textbf{0.12} & 45609.62~(0.0057) & \textbf{45607} & 4.44 \\
X-n176-k26 & \textbf{${47812^*}$} & 47937.04~(0.2615) & \textbf{47812} & 3.92 & \textbf{47812.00~(0.0000)} & \textbf{47812} & \textbf{0.41} & 47812.10~(0.0002) & \textbf{47812} & 3.76 \\
X-n181-k23 & \textbf{${25569^*}$} & 25569.54~(0.0021) & \textbf{25569} & 3.60 & \textbf{25569.00~(0.0000)} & \textbf{25569} & \textbf{1.44} & 25570.92~(0.0075) & \textbf{25569} & 4.13 \\
X-n186-k15 & \textbf{${24145^*}$} & 24158.64~(0.0565) & 24147 & 2.84 & \textbf{24145.00~(0.0000)} & \textbf{24145} & \textbf{0.48} & 24148.90~(0.0162) & \textbf{24145} & 3.44 \\
X-n190-k8 & \textbf{${16980^*}$} & 16985.72~(0.0337) & \textbf{16980} & \textbf{4.61} & 16984.16~(0.0245) & \textbf{16980} & 5.05 & \textbf{16981.18~(0.0069)} & \textbf{16980} & 4.87 \\
X-n195-k51 & \textbf{${44225^*}$} & 44262.02~(0.0837) & \textbf{44225} & 4.90 & \textbf{44225.72~(0.0016)} & \textbf{44225} & \textbf{2.62} & 44263.56~(0.0872) & \textbf{44225} & 6.33 \\
X-n200-k36 & \textbf{${58578^*}$} & 58824.86~(0.4214) & 58619 & 4.38 & \textbf{58579.44~(0.0025)} & \textbf{58578} & \textbf{3.03} & 58589.42~(0.0195) & \textbf{58578} & 5.68 \\
X-n204-k19 & \textbf{${19565^*}$}& 19566.96~(0.0100) & \textbf{19565} & 3.54 & \textbf{19565.00~(0.0000)} & \textbf{19565} & \textbf{1.42} & 19566.06~(0.0054) & \textbf{19565} & 5.74 \\
X-n209-k16 & \textbf{${30656^*}$} & 30677.96~(0.0716) & \textbf{30656} & 5.47 & \textbf{30656.70~(0.0023)} & \textbf{30656} & \textbf{3.80} & 30661.60~(0.0183) & \textbf{30656} & 5.59 \\
X-n214-k11 & \textbf{${10856^*}$} & 10877.18~(0.1951) & 10861 & \textbf{4.45} & \textbf{10863.88~(0.0726)} & \textbf{10856} & 5.91 & 10869.52~(0.1245) & 10858 & 5.01 \\
X-n219-k73 & \textbf{${117595^*}$} & 117595.14~(0.0001) & \textbf{117595} & 1.85 & 117597.42~(0.0021) & \textbf{117595} & 3.57 & \textbf{117595.00~(0.0000)} & \textbf{117595} & \textbf{1.79} \\
X-n223-k34 & \textbf{${40437^*}$} & 40511.34~(0.1838) & 40445 & 4.73 & \textbf{40438.74~(0.0043)} & \textbf{40437} & \textbf{4.15} & 40446.78~(0.0242) & \textbf{40437} & 6.19 \\
X-n228-k23 & \textbf{${25742^*}$} & 25785.5~(0.1690) & 25743 & 5.47 & \textbf{25742.74~(0.0029)} & \textbf{25742} & \textbf{2.84} & 25752.16~(0.0395) & \textbf{25742} & 5.69 \\
X-n233-k16 & \textbf{${19230^*}$} & 19294.26~(0.3342) & \textbf{19230} & \textbf{2.99} & \textbf{19230.62~(0.0032)} & \textbf{19230} & 3.45 & 19260.16~(0.1568) & \textbf{19230} & 4.39 \\
X-n237-k14 & \textbf{${27042^*}$} & 27046.34~(0.0160) & \textbf{27042} & \textbf{3.23} & \textbf{27042.16~(0.0006)} & \textbf{27042} & 3.31 & 27052.06~(0.0372) & \textbf{27042} & 5.62 \\
X-n242-k48 & \textbf{${82751^*}$} & 82898.1~(0.1778) & 82792 & 7.06 & \textbf{82808.32~(0.0693)} & \textbf{82751} & \textbf{5.29} & 82817.56~(0.0804) & \textbf{82751} & 8.56 \\
X-n247-k50 & \textbf{${37274^*}$} & 37480.2~(0.5532) & \textbf{37274} & 7.23 & 37283.46~(0.0254) & \textbf{37274} & \textbf{6.53} & \textbf{37274.32~(0.0009)} & \textbf{37274} & 8.46 \\
X-n251-k28 & \textbf{${38684^*}$} & 38785.5~(0.2624) & 38699 & \textbf{4.72} & \textbf{38693.10~(0.0235)} & \textbf{38684} & 6.45 & 38745.58~(0.1592) & \textbf{38684} & 6.73 \\
X-n256-k16 & $18839$ & 18880~(0.2176) & 18880 & \textbf{3.77} & \textbf{18841.00~(0.0106)} & \textbf{18839} & 5.61 & 18876.48~(0.1989) & \textbf{18839} & 6.12 \\
X-n261-k13 & \textbf{${26558^*}$} & 26649.02~(0.3427) & 26561 & \textbf{4.80} & \textbf{26564.94~(0.0261)} & \textbf{26558} & 7.44 & 26575.22~(0.0648) & \textbf{26558} & 6.73 \\
X-n266-k58 & \textbf{${75478^*}$} & 75766.56~(0.3823) & 75591 & \textbf{4.05} & 75609.68~(0.1745) & 75529 & 7.10 & \textbf{75543.12~(0.0863)} & \textbf{75478} & 7.22 \\
X-n270-k35 & \textbf{${35291^*}$} & 35355.46~(0.1827) & 35309 & 6.10 & \textbf{35303.00~(0.0340)} & \textbf{35291} & \textbf{3.58} & 35344.92~(0.1528) & 35303 & 7.39 \\
X-n275-k28 & \textbf{${21245^*}$} & 21251.76~(0.0318) & \textbf{21245} & 3.90 & \textbf{21245.00~(0.0000)} & \textbf{21245} & \textbf{0.60} & 21248.76~(0.0177) & \textbf{21245} & 6.14 \\
X-n280-k17 & $33503$ & 33633.2~(0.3886) & 33507 & 9.85 & \textbf{33544.48~(0.1238)} & \textbf{33503} & \textbf{8.96} & 33559.62~(0.1690) & 33505 & 9.09 \\
X-n284-k15 & \textbf{${20215^*}$} & 20269.86~(0.2714) & 20226 & 7.98 & \textbf{20249.30~(0.1697)} & 20227 & \textbf{7.55} & 20259.22~(0.2187) & \textbf{20225} & 8.46 \\
X-n289-k60 & \textbf{${95151^*}$} & 95542.14~(0.4111) & 95440 & \textbf{6.88} & 95299.00~(0.1555) & 95187 & 8.79 & \textbf{95250.44~(0.1045)} & \textbf{95157} & 8.65 \\
X-n294-k50 & $47161$ & 47264.18~(0.2188) & 47169 & \textbf{6.49} & \textbf{47191.46~(0.0646)} & \textbf{47161} & 6.80 & 47215.84~(0.1163) & 47168 & 9.59 \\
X-n298-k31 & \textbf{${34231^*}$} & 34290.3~(0.1732) & 34234 & 7.07 & \textbf{34234.14~(0.0092)} & \textbf{34231} & \textbf{5.84} & 34267.06~(0.1053) & \textbf{34231} & 9.43 \\
X-n303-k21 & $21736$ & 21809.78~(0.3394) & 21762 & 9.68 & \textbf{21746.98~(0.0505)} & 21738 & \textbf{7.77} & 21764.66~(0.1319) & \textbf{21736} & 9.07 \\
X-n308-k13 & $25859$ & 25917.76~(0.2272) & 25866 & \textbf{7.72} & \textbf{25873.44~(0.0558)} & \textbf{25861} & 7.96 & 25879.48~(0.0792) & \textbf{25861} & 9.24 \\
X-n313-k71 & $94043$ & 94327.5~(0.3025) & 94203 & 9.79 & \textbf{94107.84~(0.0689)} & \textbf{94044} & \textbf{8.55} & 94112.38~(0.0738) & \textbf{94044} & 10.25 \\
X-n317-k53 & \textbf{${78355^*}$} & 78357.26~(0.0029) & \textbf{78355} & 10.03 & 78357.56~(0.0033) & \textbf{78355} & \textbf{7.15} & \textbf{78356.66~(0.0021)} & \textbf{78355} & 8.91 \\
X-n322-k28 & \textbf{${29834^*}$} & 29937.48~(0.3469) & 29854 & 7.57 & \textbf{29850.82~(0.0564)} & \textbf{29834} & \textbf{7.22} & 29858.78~(0.0831) & 29844 & 8.97 \\
X-n327-k20 & $27532$ & 27609.92~(0.2830) & 27553 & \textbf{7.49} & \textbf{27545.76~(0.0500)} & \textbf{27532} & 8.22 & 27588.60~(0.2056) & \textbf{27532} & 9.00 \\
X-n331-k15 & \textbf{${31102^*}$} & \textbf{31103.74~(0.0056)} & 31103 & 8.13 & 31103.82~(0.0059) & \textbf{31102} & \textbf{7.63} & 31104.52~(0.0081) & 31103 & 9.73 \\
X-n336-k84 & $139111$ & 139646.86~(0.3852) & 139406 & 10.54 & 139346.04~(0.1690) & 139210 & \textbf{10.05} & \textbf{139261.64~(0.1083)} & \textbf{139140} & 12.41 \\
X-n344-k43 & $42050$ & 42204.74~(0.3680) & 42089 & 10.59 & \textbf{42075.56~(0.0608)} & \textbf{42050} & \textbf{9.00} & 42106.28~(0.1338) & 42056 & 10.02 \\
X-n351-k40 & $25896$ & 26004.68~(0.4197) & 25958 & 11.04 & 25955.08~(0.2281) & 25934 & \textbf{10.03} & \textbf{25952.10~(0.2166)} & \textbf{25905} & 12.12 \\
X-n359-k29 & $51505$ & 51624.12~(0.2313) & 51510 & \textbf{10.85} & 51661.66~(0.3042) & 51512 & 11.97 & \textbf{51553.12~(0.0934)} & \textbf{51505} & 13.07 \\
X-n367-k17 & $22814$ & 22819.4~(0.0237) & \textbf{22814} & 10.57 & \textbf{22814.02~(0.0001)} & \textbf{22814} & \textbf{4.46} & 22824.74~(0.0471) & \textbf{22814} & 11.13 \\
X-n376-k94 & \textbf{${147713^*}$} & 147721.78~(0.0059) & \textbf{147713} & \textbf{8.72} & 147718.10~(0.0035) & \textbf{147713} & 9.06 & \textbf{147714.82~(0.0012)} & \textbf{147713} & 8.80 \\
X-n384-k52 & $65928$ & 66114.66~(0.2831) & 66004 & \textbf{11.49} & 66083.22~(0.2354) & 66017 & 12.02 & \textbf{66026.56~(0.1495)} & \textbf{65938} & 11.56 \\
X-n393-k38 & \textbf{${38260^*}$} & 38307.78~(0.1249) & 38274 & 10.82 & \textbf{38260.74~(0.0019)} & \textbf{38260} & \textbf{9.37} & 38285.08~(0.0656) & \textbf{38260} & 13.05 \\
X-n401-k29 & $66154$ & 66265.74~(0.1689) & 66208 & \textbf{12.28} & 66253.22~(0.1500) & 66181 & 13.49 & \textbf{66217.14~(0.0954)} & \textbf{66170} & 15.38 \\
X-n411-k19 & $19712$ & 19785.9~(0.3749) & 19751 & 15.93 & \textbf{19724.04~(0.0611)} & \textbf{19712} & \textbf{12.75} & 19736.96~(0.1266) & 19716 & 14.48 \\
X-n420-k130 & \textbf{${107798^*}$} & 107933.82~(0.1260) & 107831 & 14.56 & 107844.12~(0.0428) & \textbf{107798} & \textbf{11.20} & \textbf{107830.74~(0.0304)} & \textbf{107798} & 16.02 \\
X-n429-k61 & $65449$ & 65589.8~(0.2151) & 65489 & 13.79 & \textbf{65501.72~(0.0806)} & \textbf{65455} & \textbf{10.88} & 65519.26~(0.1074) & 65459 & 13.19 \\
X-n439-k37 & \textbf{${36391^*}$} & 36401.04~(0.0276) & \textbf{36395} & 12.92 & \textbf{36397.66~(0.0183)} & \textbf{36395} & \textbf{10.36} & 36403.08~(0.0332) & \textbf{36395} & 13.13 \\
X-n449-k29 & $55233$ & 55417.22~(0.3335) & 55326 & 15.83 & 55403.38~(0.3085) & 55272 & \textbf{15.35} & \textbf{55307.80~(0.1354)} & \textbf{55252} & 17.41 \\
X-n459-k26 & $24139$ & 24192.04~(0.2197) & 24141 & 15.36 & \textbf{24163.60~(0.1019)} & \textbf{24139} & \textbf{12.99} & 24173.20~(0.1417) & 24141 & 13.93 \\
X-n469-k138 & \textbf{${221824^*}$} & 223116.1~(0.5825) & 222641 & \textbf{15.94} & 222211.08~(0.1745) & 221992 & 16.48 & \textbf{222026.92~(0.0915)} & \textbf{221834} & 16.28 \\
X-n480-k70 & $89449$ & 89617~(0.1878) & 89467 & 17.47 & 89546.38~(0.1089) & 89465 & \textbf{14.59} & \textbf{89467.82~(0.0210)} & \textbf{89457} & 19.28 \\
X-n491-k59 & $66483$ & 66698.44~(0.3241) & 66581 & 18.42 & 66664.16~(0.2725) & 66553 & \textbf{14.94} & \textbf{66574.48~(0.1376)} & \textbf{66487} & 20.59 \\
X-n502-k39 & $69226$ & 69244.26~(0.0264) & 69227 & 20.60 & 69249.14~(0.0334) & 69231 & \textbf{16.00} & \textbf{69233.82~(0.0113)} & \textbf{69226} & 17.95 \\
X-n513-k21 & $24201$ & 24227.38~(0.1090) & \textbf{24201} & 15.52 & \textbf{24202.62~(0.0067)} & \textbf{24201} & \textbf{10.18} & 24232.20~(0.1289) & \textbf{24201} & 17.38 \\
X-n524-k153 & \textbf{${154593^*}$} & 154908.66~(0.2042) & 154611 & 22.19 & 154765.42~(0.1115) & 154610 & \textbf{15.21} & \textbf{154631.20~(0.0247)} & \textbf{154601} & 23.39 \\
X-n536-k96 & $94828$ & 95532.1~(0.7425) & 95367 & 20.74 & 95089.34~(0.2756) & 94992 & \textbf{18.34} & \textbf{94939.66~(0.1178)} & \textbf{94878} & 22.53 \\
X-n548-k50 & \textbf{${86700^*}$} & \textbf{86738.1~(0.0439)} & \textbf{86700} & 20.19 & 86781.26~(0.0937) & 86708 & \textbf{15.97} & 86745.88~(0.0529) & \textbf{86700} & 17.68 \\
X-n561-k42 & $42717$ & 42827.5~(0.2587) & 42751 & 20.17 & \textbf{42746.98~(0.0702)} & \textbf{42719} & \textbf{16.16} & 42761.98~(0.1053) & 42723 & 20.93 \\
X-n573-k30 & $50673$ & 50802.14~(0.2548) & 50753 & 24.39 & 50805.24~(0.2610) & 50743 & 25.88 & \textbf{50735.52~(0.1234)} & \textbf{50721} & \textbf{23.41} \\
X-n586-k159 & $190316$ & 190968.62~(0.3429) & 190716 & 22.86 & 190589.98~(0.1440) & 190447 & \textbf{19.20} & \textbf{190399.76~(0.0440)} & \textbf{190318} & 24.87 \\
X-n599-k92 & $108451$ & 108744.02~(0.2702) & 108645 & 23.11 & 108697.96~(0.2277) & 108551 & \textbf{19.44} & \textbf{108589.12~(0.1274)} & \textbf{108484} & 23.89 \\
X-n613-k62 & $59535$ & 59714.54~(0.3016) & 59585 & 24.10 & 59697.26~(0.2725) & 59585 & \textbf{16.31} & \textbf{59611.80~(0.1290)} & \textbf{59545} & 24.40 \\
X-n627-k43 & $62164$ & 62279.3~(0.1855) & 62188 & 26.06 & 62380.94~(0.3490) & 62234 & 27.27 & \textbf{62222.80~(0.0946)} & \textbf{62179} & \textbf{25.93} \\
X-n641-k35 & $63682$ & 63822.38~(0.2204) & \textbf{63734} & 25.73 & 63893.24~(0.3317) & 63750 & 24.07 & \textbf{63777.22~(0.1495)} & \textbf{63734} & \textbf{23.41} \\
X-n655-k131 & \textbf{${106780^*}$} & 106807.04~(0.0253) & \textbf{106780} & 22.48 & 106813.18~(0.0311) & 106785 & \textbf{19.06} & \textbf{106793.58~(0.0127)} & \textbf{106780} & 20.87 \\
X-n670-k130 & $146332$ & 147001.6~(0.4576) & 146839 & 31.70 & 146913.28~(0.3972) & 146517 & \textbf{23.27} & \textbf{146831.38~(0.3413)} & \textbf{146419} & 31.21 \\
X-n685-k75 & $68205$ & 68466.92~(0.3840) & 68310 & 27.56 & 68394.94~(0.2785) & 68293 & \textbf{22.45} & \textbf{68304.10~(0.1453)} & \textbf{68228} & 28.92 \\
X-n701-k44 & $81919$ & 82116.36~(0.2409) & 81965 & 30.64 & 82267.08~(0.4249) & 82040 & \textbf{27.18} & \textbf{81985.30~(0.0809)} & \textbf{81923} & 30.39 \\
X-n716-k35 & $43356$ & 43491.7~(0.3130) & 43430 & 31.45 & 43513.22~(0.3626) & 43451 & \textbf{25.91} & \textbf{43401.94~(0.1060)} & \textbf{43338} & 30.84 \\
X-n733-k159 & $136187$ & 136428.9~(0.1776) & 136310 & 31.17 & 136454.62~(0.1965) & 136322 & \textbf{25.24} & \textbf{136301.84~(0.0843)} & \textbf{136231} & 30.47 \\
X-n749-k98 & $77269$ & 77583.86~(0.4075) & 77478 & 30.10 & 77697.62~(0.5547) & 77527 & \textbf{29.63} & \textbf{77420.64~(0.1962)} & \textbf{77337} & 33.77 \\
X-n766-k71 & $114416$ & 114757.88~(0.2988) & 114642 & 34.54 & 114779.78~(0.3179) & 114648 & \textbf{30.48} & \textbf{114473.72~(0.0504)} & \textbf{114423} & 35.17 \\
X-n783-k48 & $72381$ & 72601.54~(0.3047) & 72471 & \textbf{30.58} & 72832.96~(0.6244) & 72641 & 32.90 & \textbf{72481.30~(0.1386)} & \textbf{72400} & 33.16 \\
X-n801-k40 & $73305$ & 73400.48~(0.1303) & 73320 & \textbf{31.22} & 73487.62~(0.2491) & 73345 & 32.07 & \textbf{73366.04~(0.0833)} & \textbf{73308} & 31.26 \\
X-n819-k171 & $158121$ & 158925.12~(0.5085) & 158740 & 36.11 & 158498.82~(0.2389) & 158320 & \textbf{29.36} & \textbf{158254.84~(0.0846)} & \textbf{158187} & 36.28 \\
X-n837-k142 & $193734$ & 194225.32~(0.2536) & 194059 & \textbf{33.78} & 194278.88~(0.2813) & 194045 & 38.02 & \textbf{193832.18~(0.0507)} & \textbf{193748} & 36.95 \\
X-n856-k95 & \textbf{${88965^*}$} & 89030.32~(0.0734) & 88973 & 34.97 & 89038.14~(0.0822) & 88983 & \textbf{24.47} & \textbf{89019.78~(0.0616)} & \textbf{88966} & 31.93 \\
X-n876-k59 & $99299$ & 99514.72~(0.2172) & 99407 & \textbf{39.52} & 99731.38~(0.4354) & 99501 & 42.90 & \textbf{99427.90~(0.1298)} & \textbf{99352} & 41.12 \\
X-n895-k37 & $53848$ & 54037.02~(0.3510) & 53901 & 35.62 & 54134.80~(0.5326) & 53975 & \textbf{30.30} & \textbf{53956.02~(0.2006)} & \textbf{53881} & 34.35 \\
X-n916-k207 & $329178$ & 330184.72~(0.3058) & 329942 & \textbf{41.53} & 329966.70~(0.2396) & 329719 & 44.54 & \textbf{329316.82~(0.0422)} & \textbf{329205} & 42.51 \\
X-n936-k151 & $132715$ & 133430.2~(0.5389) & 133047 & 45.21 & 133416.98~(0.5289) & 133118 & \textbf{30.29} & \textbf{132991.78~(0.2086)} & \textbf{132909} & 43.08 \\
X-n957-k87 & $85464$ & 85527.7~(0.0745) & 85474 & 37.83 & 85564.24~(0.1173) & 85494 & \textbf{30.96} & \textbf{85511.52~(0.0556)} & \textbf{85473} & 35.06 \\
X-n979-k58 & $118954$ & 119289.54~(0.2821) & 119195 & \textbf{43.91} & 119368.28~(0.3483) & 119202 & 48.55 & \textbf{119010.64~(0.0476)} & \textbf{118951} & 46.57 \\
X-n1001-k43 & $72355$ & 72547.88~(0.2666) & 72441 & 44.88 & 72757.08~(0.5557) & 72598 & 46.65 & \textbf{72455.22~(0.1385)} & \textbf{72375} & \textbf{44.50} \\
Average & & ~(0.2049) & 14.11 & & ~(0.1268) & & 12.49 & ~(0.0782) & & 14.87 \\* \bottomrule
\end{longtable}
\end{landscape}

\section{Parameter analysis}\label{appendix2}

\textcolor{black}{This section shows a summarized analysis of four parameters of AILS-II: $d_{max}$, $d_{min}$, $\gamma$ and $\varphi$. In this experiment, we evaluated these parameters considering the following sets of values:}

\begin{itemize}
    \item \textcolor{black}{$d_{max}\in\{10,20,30,40,50,60,70,80\}$};
     \item \textcolor{black}{$d_{min}\in \{5,10,15,20,25,30\}$};
      \item \textcolor{black}{$\gamma\in\{10,20,30,40,50,60\}$};
       \item \textcolor{black}{$\varphi\in\{10,20,30,40,50,60,70,80\}$}.
\end{itemize}

\textcolor{black}{For this sensitivity analysis, we employed the  values of $d_{max}=30, d_{min}=15,\gamma=30$ and $\varphi=40$ when varying one of the parameters. Figure~\ref{fig:variandoParam} shows the average gaps on the large-scale instances.}
\begin{figure}[h!t]
\center
\includegraphics[]{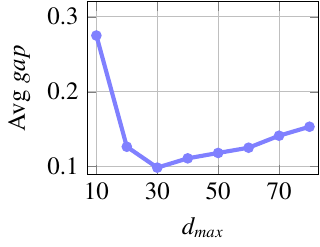}
\includegraphics[]{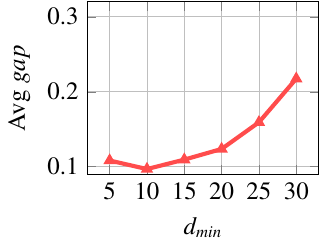}\\
\includegraphics[]{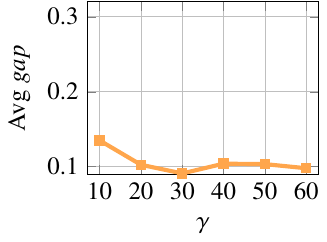}
\includegraphics[]{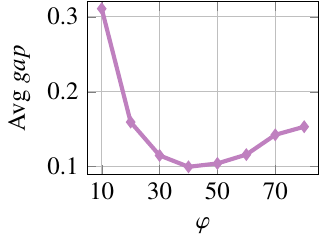}
\caption{\textcolor{black}{Average results achieved by varying the parameters of AILS-II considering the instances of \cite{Arnold2019}.} }
\label{fig:variandoParam}
\end{figure}

\textcolor{black}{It is noteworthy that, except for $\gamma$, the other parameters were sensitive to variations, with abrupt changes in the average gap. The final values reported in the experiment were achieved when $d_{max}=30, d_{min}=15,\gamma=30$ and $\varphi=40$. Lower values of $\varphi$ were responsible for much poorer results. This parameter refers to the number of investigated neighbors of a given node during the search process. This particular result demonstrates the importance of limiting the search to achieve a better performance.}

 \end{APPENDIX}

\end{document}